\documentclass[10pt]{amsart}
\author{Chiung-ju Liu}
\keywords{Asymptotic expansion; Riemann surfaces}
\title{The asymptotic Tian-Yau-Zelditch expansion on Riemann surfaces with Constant Curvature}
\email{cjliu4@ntu.edu.tw}
\address{Taida Institute of Mathematical Science,
National Taiwan University, Taiwan 106}
\date{June 30, 2008}
\thanks{Research supported by the grant of NSC A95-0918.}
\usepackage{amsfonts}
\usepackage{amsmath,amsthm,amscd, amssymb}
\usepackage{amsthm}
\usepackage{setspace}
\usepackage{amsrefs}
\begin{document}
 \footnote{\textit{Subject Classification.} Primary
53C55; Secondary 32Q05\nonumber}
\begin{abstract}
 Let $M$ be a regular Riemann surface with a
 metric which has constant scalar curvature $\rho$. We give
the asymptotic expansion of the sum of the square norm of the
sections of the pluricanonical bundles $K_{M}^{m}$. That is,
\[\sum_{i=0}^{d_{m}-1}\|S_{i}(x_{0})\|_{h_{m}}^{2}
\sim m(1+\frac{\rho}{2 m})+O\left(e^{-\frac{(\log
m)^{2}}{8}}\right),\] where $\{S_{0},\cdots,S_{d_{m}-1}\}$ is an
orthonormal basis for $H^{0}(M, K_{M}^{m})$ for sufficiently large
$m$.
\end{abstract}
\newtheorem{thm}{Theorem}[section]
\maketitle \section{Introduction} Let $M$ be an $n$-dimensional
compact complex K\"ahler manifold with an ample line bundle $L$
over $M$. Let $g$ be the K\"ahler metric on $M$ corresponding to
the K\"ahler form $\omega_{g}=Ric(h)$ for some positive Hermitian
$h$ metric on $L$. Such a K\"ahler metric $g$ is called a
polarized K\"ahler metric. The metric $h$ induces a Hermitian
metric $h_{m}$ on $L^{m}$ for all positive integers $m$. Let
$\{S_{0},\cdots,S_{d_{m}-1}\}$ be an orthonormal basis of the
space $H^{0}(M,L^{m})$ with respect to the inner product
\numberwithin{equation}{section}
\begin{equation}(S,T)=\int_{M}\langle
S(x),T(x)\rangle_{h_{m}}dV_{g},\label{2}\end{equation} where
$d_{m}=\mathrm{dim}H^{0}(M,L^{m})$ and
$dV_{g}=\frac{\omega_{g}^{n}}{n!}$ is the volume form of $g$. The
quantity
\begin{equation}\sum_{i=0}^{d_{m}-1}\|S_{i}(x)\|^{2}_{h_{m}}\label{1}\end{equation} is
related to the existence of K\"ahler-Einstein metrics and
stability of complex manifolds. A lot of work has been done for
\eqref{1} on compact complex K\"ahler manifolds. Tian~\cite{TG5}
applied H\"omander's $L^{2}$-estimate to produce peak sections and
proved the $C^{2}$ convergence of the Bergman metrics. Later,
Ruan~\cite{RWD} proved the $C^{\infty}$ convergence. About the
same time, Zelditch~\cite{ZS} and Catlin~\cite{CD} separately
generalized the theorem of Tian by showing there is an asymptotic
expansion
\begin{equation}
\sum_{i=0}^{d_{m}-1}\|S_{i}(x)\|^{2}_{h_{m}}\sim
a_{0}(x)m^{n}+a_{1}(x)m^{n-1}+a_{2}(x)m^{n-2}+\cdots
\label{1.1}\end{equation} for certain smooth coefficients $a_{j}(x)$
with $a_{0}=1$.
 In~\cite{LZ4}, Lu proved that each coefficient
$a_{j}(x)$ is a polynomial of the curvature and its covariant
derivatives. In particular, $a_{1}=\frac{\rho}{2}$, where $\rho$
is the scalar curvature of $M$. These polynomials can be found by
finitely many steps of algebraic operations. Recently, Song
\cite{S} generalized Zelditch's theorem on orbifolds of finite
isolated singularities. The information on the singularities can
be found in the expansion.

On the Riemann surfaces with bounded curvature, Lu~\cite{LZ3}
proved that there is a lower bound for \eqref{1}. Later, the
result of Lu and Tian~\cite{L-T} implies that on the Riemann
surfaces with constant scalar curvature $\rho$, the asymptotic
expansion \eqref{1.1} is given by
\[\sum_{i=0}^{d_{m}-1}\|S_{i}(x_{0})\|_{h_{m}}^{2}
\sim m(1+\frac{\rho}{2m})+O\left(\frac{1}{m^{p}}\right)\] for any
$p>0$. In the current paper, we obtain a more precise result for
\eqref{1.1}.

\begin{thm}\label{Liu}
Let $M$ be a regular compact Riemann surface and $K_{M}$ be the
canonical line bundle endowed with a Hermitian metric $h$ such
that the curvature $Ric(h)$ of $h$ defines a K\"ahler metric $g$
on $M$. Suppose that this metric $g$ has constant scalar curvature
$\rho$. Then there is a complete asymptotic expansion:
\[\sum_{i=0}^{d_{m}-1}\|S_{i}(x_{0})\|_{h_{m}}^{2}
\sim m(1+\frac{\rho}{2m})+O\left(e^{-\frac{(\log
m)^{2}}{8}}\right),\] where $\{S_{0},\cdots,S_{d_{m}-1}\}$ is an
orthonormal basis for $H^{0}(M, K_{M}^{m})$ for some
$m>\max\{e^{20\sqrt{5}}+2|\rho|,|\rho|^{4/3},\frac{1}{\delta},\sqrt{\frac{2}{|\rho|}}\}$,
where $\delta$ is the injective radius at $x_0$.
\end{thm}
Our result indicates that the asymptotic expansion \eqref{1.1} is
in exponential decay.\\
Engli\u s~\cite{EM} has an asymptotically expansion for the
Berezin transformation on any planar domain of hyperbolic type. He
also showed that Berezin kernel~\cite{B} has
\[\tilde{B}(\eta,\eta)=m\left(1+O(1)\rho_{0}(0)^{\frac{\pi m-3}{2}}\right),\]
where $\rho_{0}(0)$ is a positive constant.

\paragraph{\textbf{Acknowledgements}}
The author thanks her advisor Z. Lu for his suggestions and help
during the preparation of this paper. She also thanks R. Chiang
for her encouragement.
\section{General set up}

Let $M$ be an $n$-dimensional compact complex K\"ahler manifold with
a polarized line bundle $(L,h)\rightarrow M$. Choose the
$K$-coordinates $(z_{1},\cdots,z_{n})$ on an open neighborhood $U$
of a fixed point $x_{0}\in M$. Then the K\"ahler form
\[\omega_{g}=\frac{\sqrt{-1}}{2\pi}\sum_{\alpha,\beta=1}^{n}g_{\alpha\bar{\beta}}dz_{\alpha}\wedge
d\bar{z}_{\beta}\] satisfies
\begin{equation}
g_{\alpha\bar{\beta}}(x_{0})=\delta_{\alpha\bar{\beta}}, \qquad
\frac{\partial^{p_{1}+\cdots+p_{n}}g_{\alpha\bar{\beta}}}{\partial
z_{1}^{p_{1}}\cdots\partial z_{n}^{p_{n}}}(x_{0})=0,\label{2.1}
\end{equation}
for $\alpha,\beta=1,\cdots,n$ and any nonnegative integers
$p_{1},\cdots, p_{n}$ with $p_{1}+\cdots+p_{n}\neq 0$.

We also choose a local holomorphic frame $e_{L}$ of the line bundle
$L$ at $x_{0}$ such that $a$ is the local representation function of
the Hermitian metric $h$. That is,
\[Ric(h)=-\frac{\sqrt{-1}}{2\pi}\partial\bar{\partial}\log
a.\] Under the $K$-coordinate, the function $a$ has the properties
\begin{equation}
a(x_{0})=1,\qquad \frac{\partial^{p_{1}+\cdots+p_{n}}}{\partial
z_{1}^{p_{1}}\cdots\partial z_{n}^{p_{n}}}(a)(x_{0})=0\label{2.2}
\end{equation}
for any nonnegative integers $p_{1},\cdots, p_{n}$ with
$p_{1}+\cdots+p_{n}\neq 0$.

Let $\{S_{0},\cdots,S_{d_{m}-1}\}$ be a basis of $H^{0}(M,L^{m})$.
Assume that at the point $x_{0}\in M$,
\[S_{0}(x_{0})\neq 0, \qquad S_{i}(x_{0})=0, \qquad
i=1,\cdots,d_{m}-1.\] If the set $\{S_{0},\cdots,S_{d_{m}-1}\}$ is
not an orthonormal basis, we may do the following: \\Let the
metric matrix
\[F_{ij}=(S_{i},S_{j}), \qquad i,j=0,\cdots, d_{m}-1\] with respect
to the inner product \eqref{2}. By definition, $(F_{ij})$ is a
positive definite Hermitian matrix. We can find a $d_{m}\times
d_{m}$ matrix $G_{ij}$ such that
\[F_{ij}=\sum_{k=0}^{d_{m}-1}G_{ik}\overline{G_{jk}}.\] Let $(H_{ij})$
be the inverse of $(G_{ij})$. Then
$\{\sum_{j=0}^{d_{m}-1}H_{ij}S_{j}\}$ forms an orthonormal basis of
$H^{0}(M,L^{m})$. The left hand side of \eqref{1} is equal to
\begin{equation}
\sum_{i=0}^{d_{m}-1}\|\sum_{j=0}^{d_{m}-1}H_{ij}S_{j}(x_{0})\|^{2}_{h_{m}}=
\sum_{i=0}^{d_{m}-1}|H_{i0}|^{2}\|S_{0}(x_{0})\|^{2}_{h_{m}}.\label{2.3}
\end{equation}
Let $(I_{ij})$ be the inverse matrix of $(F_{ij})$. Denote that
\begin{equation}\sum_{i=0}^{d_{m}-1}|H_{i0}|^{2}=I_{00}.\label{2.4}
\end{equation}
In order to compute \eqref{2.4}, we need a suitable choice of the
basis $\{S_{0},\cdots, S_{d_{m}-1}\}$. We select some of Tian's
peak sections in our basis. The following lemma an improved
version of Tian's result~\cite[Lemma 1.2]{TG5}, which is done by
Lu and Tian.

Let $\mathbb{Z}_{+}^{n}$ be the set of $n$-tuple integers
$P=(p_{1},\cdots,p_{n})$ such that each $p_{i}$ is a nonnegative
integer for $i=1,\cdots,n$. For $P\in\mathbb{Z}_{+}^{n}$, we denote
that $z^{P}=z_{1}^{p_{1}}\cdots z_{n}^{p_{n}}$ and $|P|=p_{1}+\cdots
p_{n}$.
\newtheorem{lemma}{Lemma}[section]
\begin{lemma}[Tian]\label{Tian}
Suppose $Ric(g)\geq -K\omega_{g}$, where $K>0$ is a constant. For
$P\in\mathbb{Z}_{+}^{n}$ and an integer $p'>|P|$, let $m$ be an
integer such that \[m>\max\{e^{20\sqrt{n+2p'}}+2K,
e^{8(p'-1+n)}\}.\] Then there is a holomorphic section $S_{P,m}\in
H^{0}(M,L^{m})$, satisfying
\begin{equation}
|\int_{M}\|S_{P,m}\|_{h_{m}}^{2}dV_{g}-1|\leq C e^{-\frac{1}{8}(\log
m)^{2}}.\label{2.5}
\end{equation}
Moreover, $S_{P,m}$ can be decomposed as
\[S_{P,m}=\tilde{S}_{P,m}-u_{P,m}\]
such that
\begin{equation}
\tilde{S}_{P,m}(x)=\lambda_{P}\eta\left(\frac{m|z|^{2}}{(\log
m)^{2}}\right)z^{P}e_{L}^{m}=
\begin{cases}\lambda_{P}z^{P}e_{L}^{m}&\text{$x\in\{|z|\leq \frac{\log m}{\sqrt{2m}}\}$},\\
0&\text{$x\in M\setminus \{|z|\leq \frac{\log
m}{\sqrt{m}}\}$},\label{2.7}
\end{cases}
\end{equation}
and \begin{equation}\int_{M}\|u_{P,m}\|^{2}_{h_{m}}dV_{g}\leq C
e^{-\frac{1}{4}(\log m)^{2}},\label{2.9}
\end{equation}
where $\eta$ is a smoothly cut-off function
\begin{equation*}\eta(t)=\begin{cases}1&\text{for $0\leq t \leq \frac{1}{2}$},\\0&\text{for $t\geq 1$}.\end{cases}
\end{equation*} satisfying $0\leq -\eta'(t)\leq 4$ and $|\eta''(t)|\leq 8$ and
\begin{equation}
\lambda_{P}^{-2}=\int_{|z|\leq\frac{\log
m}{\sqrt{m}}}|z^{P}|^{2}a^{m}dV_{g}.\label{2.6}
\end{equation}
\end{lemma}
\begin{proof}  Define the weight function
\[\Psi(z)=(n+2p')\eta\left(\frac{m|z|^{2}}{(\log m)^{2}}\right)\log\left(\frac{m|z|^{2}}{(\log m)^{2}}\right).\]
A straightforward computation gives
\begin{equation}
\sqrt{-1}\partial\bar{\partial}\Psi\geq -\frac{100m(n+2p')}{(\log
m)^{2}}\omega_{g}.\label{2.10}
\end{equation}
By using \eqref{2.10}, we can verify that
\[\langle \partial\bar{\partial}\Psi+\frac{2\pi}{\sqrt{-1}}(Ric(h^{m})+Ric(g)),v\wedge \bar{v}\rangle_{g}\geq\frac{1}{4}m\|v\|_{g}^{2}.\]
For $P\in \mathbb{Z}_{+}^{n}$, consider the $1$-form
\[w_{P}=\bar{\partial}\eta(\frac{m|z|^{2}}{(\log
m)^{2}})z^{P}e_{L}^{m}.\] Since $w_{P}\equiv 0$ in a neighborhood
of $x_{0}$, we have
\[\int_{M}\|w_{P}\|_{h_{m}}^{2}e^{-\Psi}dV_{g}<+\infty.\]
By~\cite[Prop. 2.1]{TG5}, there exists a smooth $L^{m}$-valued
section $u_{P}$ such that $\bar{\partial}u_{P}=w_{P}$ and
\begin{equation}
\int_{M}\|u_{P}\|^{2}_{h_{m}}e^{-\Psi}dV_{g}
\leq\frac{4}{m}\int_{M}\|w_{P}\|_{h_{m}}^{2}e^{-\Psi}dV_{g}<\infty.\label{2.11}
\end{equation}
By direct computation, we get
\[\int_{M}\|u_{P}\|_{h_{m}}^{2}e^{-\Psi}dV_{g}\leq
\frac{C(\log m)^{2(p-1)}}{m^{p}}\int_{\frac{\log
m}{\sqrt{2m}}\leq|z|\leq \frac{\log m}{\sqrt{m}}}a^{m}dV_{0}.\]
Under the $K$-coordinate, we have
\[a^{m}=e^{m\log a}=e^{m(-|z|^{2}+O(|z|^{4}))}.\]
Hence we get
\[\int_{M}\|u_{P}\|_{h_{m}}^{2}e^{-\Psi}dV_{g}\leq \frac{C_{1}(\log
m)^{2(p-1+n)}}{m^{p+n}}e^{-\frac{1}{2}(\log m)^{2}}\] for some
constant $C_{1}$.
 Let
$\tilde{S}_{P,m}=\lambda_{P}\eta(\frac{m|z|^{2}}{(\log
m)^{2}})z^{P}e_{L}^{m}$ and $u_{P,m}=\lambda_{P}u_{P}$. Use a
result in~\cite{LZ4} \[\lambda_{P}^{2}\leq C_{2}m^{n+|P|}\] for
some constant $C_{2}$. Then we have
\[\int_{M}\|u_{P,m}\|_{h_{m}}^{2}dV_{g}\leq C(\log
m)^{2(|P|-1+n)}e^{-\frac{1}{2}(\log m)^{2}}.\] Choosing
$m>e^{8(p'-1+n)}$, we obtain
\[\int_{M}\|u_{P,m}\|_{h_{m}}^{2}dV_{g}\leq Ce^{-\frac{1}{4}(\log m)^{2}}.\]
\end{proof}
\section{Proof of Theorem \ref{Liu}}

\begin{proof}
Let $M$ be a smooth compact Riemann surface with a metric $g$ that
has constant scalar curvature. Let $x_0$ be a fixed point. Let
\[U=\{x:{\rm dist}(x,x_0)< \delta\},\] where $\delta$ is the
injective radius at $x_0$. It is well known that on a Riemann
surface there is an isothermal coordinate at each point on $U$. We
may assume that there is a holomorphic function $z$ on $U$ and it
defines the conformal structure on $U$. That is,
\[ds^2=g dz d\bar{z}\] and $g>0$.
The metric $g$ satisfies
\begin{equation} \triangle\log g=-\rho,
\quad g(x_{0})=1, \quad\text{and}\quad \frac{\partial g}{\partial
z}(x_{0})=0, \label{2.12}\end{equation} where
\[\triangle=g^{-1}\frac{\partial^{2}}{\partial z\partial\bar{z}}\]
is the complex Laplace of $M$. Since the metric $g$ has conformal
structure, it is rotationally symmetric. We can write \eqref{2.12}
in polar coordinates $(r,\theta)$:
\begin{equation} \frac{\partial^{2}g}{\partial
r^{2}}+\frac{1}{r}\frac{\partial g}{\partial
r}-\frac{1}{g}(\frac{\partial g}{\partial r})^{2}=-4\rho g^{2},
\qquad g(0,\theta)=1,\quad \frac{\partial g}{\partial
r}(0,\theta)=0,\label{2.13}
\end{equation} where $z=re^{i\theta}$, and $|z|^{2}=r^{2}$. There
exists a solution
\begin{equation}g=\frac{1}{(1+\frac{\rho}{2}|z|^{2})^{2}}\label{2.14}
\end{equation}
to \eqref{2.13} for $|z|<\sqrt{-\frac{2}{\rho}}$ if $\rho<0$.
Suppose that there exists another solution $g_{1}$ to
\eqref{2.13}. We have
\begin{equation*}
\triangle\log \left(g_{1}/g\right)=0\quad\text{and}\quad
g_1(x_0)=1.
\end{equation*}
For $\rho<0$, let $r_0<\sqrt{-\frac{2}{\rho}}$. Since $g$ and
$g_1$ are rotationally symmetric, they remain constant on
$|z|=r_0$. The harmonic function $\log(g/g_1)$ is a constant on
$|z|\leq r_0$ by Maximum Principle. By definition, we have
$g(x_{0})=g_1 (x_{0})=1$. Therefore, the solution in \eqref{2.14}
is unique around $x_{0}$. By the same reason, $g=g_1$ on $\{{\rm
dist}(x,x_0)\leq \delta_1\}$ for some $\delta_1 <\delta$ for
$\rho\leq 0$.

Let $a$ be the local representation of the metric $h$ on $K_{M}$
such that
\[-\frac{\sqrt{-1}}{2\pi}\partial\bar{\partial}\log a=\omega_{g}.\]
If we normalize $a$ and $a$ satisfies
\begin{equation}
\triangle\log a=-1, \quad a(x_{0})=1, \quad \frac{\partial
a}{\partial z}(x_{0})=0.\label{a}
\end{equation}
Since
\[-\frac{\partial^{2}}{\partial z\partial\bar{z}} \log a=g,\]
$\log a$ is also rotationally symmetric. Since
\begin{equation}a=\left\{%
\begin{array}{ll}
    \left(1+\frac{\rho}{2}|z|^{2}\right)^{-\frac{2}{\rho}}, & \hbox{if $\rho\neq 0$;} \\
    e^{-|z|^{2}}, & \hbox{if $\rho=0$.} \\
\end{array}%
\right.\label{3.5}
\end{equation} satisfies \eqref{a}, the local
uniqueness is due to the same reason.

We need to choose sufficient large $m$ such that $\frac{\log
m}{\sqrt{m}}<\min\{\delta, \sqrt{\frac{2}{|\rho|}}\}$. With these
particular solutions of $g$ and $a$, we further compute
\begin{eqnarray}\lambda_{0}^{-2}&=&\int_{|z|\leq\frac{\log
m}{\sqrt{m}}}a^{m}g\frac{\sqrt{-1}}{2\pi}dz\wedge
d\bar{z}\nonumber\\
&=&2\int_{o}^{\frac{\log
m}{\sqrt{m}}}(1+\frac{\rho}{2}r^{2})^{-\frac{2m}{\rho}-2}rdr\nonumber\\
&=&\frac{1}{m+\frac{\rho}{2}}\left(1-\big(1+\frac{\rho}{2}\frac{(\log
m)^{2}}{m}\big)^{-1-\frac{2m}{\rho}}\right)\qquad\text{for
$\rho\neq 0$}.
\end{eqnarray}
For $m>\max\{|\rho|^{4/3},10\}$, we have
$\big|\frac{\rho}{2}\frac{(\log m)^{2}}{m}\big|<1/2$. For $\rho\neq
0$, this gives
\[
\big(1+\frac{\rho}{2}\frac{(\log
m)^{2}}{m}\big)^{-1-\frac{2m}{\rho}}\leq
2e^{-\frac{2m}{\rho}\big(\frac{\rho}{2}\frac{(\log
m)^{2}}{m}-\frac{1}{2}(\frac{\rho}{2}\frac{(\log
m)^{2}}{m})^{2}+\cdots\big)}\leq Ce^{-(\log m)^{2}}.\] For
$\rho=0$, we have
\[\lambda_{0}^{-2}=\int_{|z|\leq\frac{\log
m}{\sqrt{m}}}e^{-m|z|^{2}}\frac{\sqrt{-1}}{2\pi}dz\wedge
d\bar{z}=\frac{1}{m}(1+ O(e^{-(\log m)^{2}})).\] Hence we obtain
\begin{equation}
\lambda_{0}^{-2}=
    \frac{1}{m+\frac{\rho}{2}}\left(1+O\big(e^{-(\log
m)^{2}}\big)\right). \label{2.20}\end{equation}

From the properties of $g$ and $a$, the isothermal coordinate
$(U,z)$ is a $K$-coordinate. According to Lemma \ref{Tian}, we may
choose two peak sections
\begin{eqnarray*}
S_{0,m}&=&\lambda_{0}(\eta(\frac{m|z|^{2}}{(\log
m)^{2}})(dz)^{m}-u_{0})\\
S_{1,m}&=&\lambda_{1}(\eta(\frac{m|z|^{2}}{(\log
m)^{2}})z(dz)^{m}-u_{1})\end{eqnarray*} in $H^{0}(M,K_{M}^{m})$
for some $m>e^{20\sqrt{1+4}}+2|\rho|$. Obviously, we have
$S_{0,m}(x_{0})\neq 0$ and $S_{1,m}(x_{0})=0$. Let the subspace
\[V=\{S\in H^{0}(M,K_{M}^{m})|S(x_{0})=0, DS(x_{0})=0\},\]
where $D$ is a covariant derivative on $K_{M}^{m}$. Let
$T_{1},\cdots,T_{d_{m}-2}$ be an orthonormal basis of $V$. Let
\begin{align}
S_{i}=\begin{cases}
S_{i,m} &\text{if $i=0,1$}\\
T_{i-1}&\text{if $2\leq i\leq d_{m}-1$}
\end{cases}.\label{2.15}
\end{align}
Then $\{S_{i}\}_{i=0}^{d_{m}-1}$ forms a basis for
$H^{0}(M,K_{M}^{m})$. Locally, each $T_{i}$ has the form
$f_{i}(dz)^{m}$ for some holomorphic function $f_{i}$ defined in
$U$. The holomorphic function $f_{i}$ has Taylor expansion as
$f_{i}=\sum_{\alpha=2}^{\infty}b_{i\alpha}z^{\alpha}$ in $U$,
since $T_{i}(x_{0})=0$ and $DT_{i}(x_{0})=0$ for $1\leq i\leq
d_{m}-2$
\begin{lemma}\label{matrix}
Let $\{S_{i}\}_{i=0}^{d_{m}-1}$ be the basis of $
H^{0}(M,K_{M}^{m})$, defined in \eqref{2.15}. For
$m>e^{20\sqrt{5}}+2|\rho|$, the Hermitian matrix
\[(S_{i},S_{j})=\int_{M}\langle
S_{i}(x),S_{j}(x)\rangle_{h_{m}}dV_{g}\] is given by
\begin{eqnarray*}(S_{0},S_{0})&=&1+O\left(e^{-\frac{(\log
m)^{2}}{8}}\right),\\
(S_{0},S_{1})&=&O\left(e^{-\frac{(\log m)^{2}}{8}}\right),\\
(S_{1},S_{1})&=&1+O\left(e^{-\frac{(\log m)^{2}}{8}}\right),\\
(S_{0},S_{i})&=&O\left(e^{-\frac{(\log
m)^{2}}{8}}\right),\\
(S_{1},S_{i})&=&O\left(e^{-\frac{(\log
m)^{2}}{8}}\right),\\
(S_{i},S_{j})&=&\delta_{ij}\end{eqnarray*} for
$i,j=2,\cdots,d_{m}-1$.
\end{lemma}
\begin{proof}
By definition, we have $(S_{i},S_{j})=\delta_{ij}$ for $2\leq
i,j\leq d_{m}-1$. The inner product of $(S_{i},S_{i})$ for $0\leq
i\leq 1$ is directly from Lemma \ref{Tian}. Since $a^{m}g$ is
rotationally symmetric, we have
\[\int_{|z|\leq\frac{\log
m}{\sqrt{m}}}\bar{z}^{\alpha}a^{m}gdV_{0}=0\qquad\text{for
arbitrary positive integer $\alpha$}.\] Then we get
\begin{eqnarray*}
(S_{0},S_{1})&=&(\tilde{S}_{0},\tilde{S}_{1})+(\lambda_{0}u_{0},\tilde{S}_{1})+(\tilde{S}_{0},\lambda_{1}u_{1})+(u_{0},u_{1})\\&=&
O\left(e^{-\frac{(\log m)^{2}}{8}}\right).\end{eqnarray*} Consider
\begin{eqnarray*}
(S_{0},S_{i})&=&
\int_{M}\langle\lambda_{0}(\eta(\frac{m|z|^{2}}{(\log
m)^{2}})(dz)^{m}-u_{0}),f_{i-1}(dz)^{m}\rangle_{h_{m}}dV_{g}\\
&\leq&\lambda_{0}\int_{|z|\leq\frac{\log
m}{\sqrt{m}}}\sum_{\alpha=2}^{\infty}b_{(i-1)\alpha}\bar{z}^{\alpha}a^{m}gdV_{0}+\lambda_{0}\|u_{0}\|\cdot\|S_{i}\|.
\end{eqnarray*}
Thus we have
\[(S_{0},S_{i})=O\left(e^{-\frac{(\log m)^{2}}{8}}\right)\qquad\text{for $2\leq i\leq
d_{m}-1$}.\] Similarly, consider \[(S_{1},S_{j})\leq
\lambda_{0}\int_{|z|\leq\frac{\log
m}{\sqrt{m}}}\sum_{\alpha=2}^{\infty}b_{(i-1)\alpha}z\bar{z}^{\alpha}a^{m}gdV_{0}+\lambda_{1}\|u_{1}\|\cdot\|S_{i}\|\qquad\text{
for $2\leq i\leq d_{m}-1$}.\]Since $a^{m}g$ is rotationally
symmetric, $\int_{|z|\leq\frac{\log
m}{\sqrt{m}}}z\bar{z}^{\alpha}a^{m}gdV_{0}=0$ for $\alpha\geq 2$.
Hence we obtain
\[(S_{1},S_{i})=O\left(e^{-\frac{(\log m)^{2}}{8}}\right).\]
\end{proof}
 According to \cite[Definition 3.1]{LZ4}, the metric matrix
$(F_{ij})$ can be represented by the block matrices
\begin{equation}(F_{ij})=\left(\begin{array}{ccc}
(S_{0},S_{0})&(S_{0},S_{1})&M_{13}\\
(S_{1},S_{0})&(S_{1},S_{1})&M_{23}\\
M_{31}&M_{32}&E
\end{array}\right),\label{2.16}
\end{equation}
where $M_{13}=((S_{0},S_{2}),\cdots,(S_{0},S_{d_{m}-1}))$,
$M_{23}=((S_{1},S_{2}),\cdots,(S_{1},S_{d_{m}-1}))$,
$M_{31}=\overline{M_{13}^{T}}$, $M_{32}=\overline{M_{23}^{T}}$,
and $E$ is a $(d_{m}-2)\times(d_{m}-2)$ identity matrix. By
using~\cite[Lemma 3.1]{LZ4}, we obtain
\begin{equation}I_{00}=\frac{1}{(S_{0},S_{0})}+(\frac{1}{(S_{0},S_{0})})^{2}\left(\begin{array}{cc}
(S_{0},S_{1})&M_{13}\end{array}\right)\tilde{M}^{-1}\left(\begin{array}{c}(S_{1},S_{0})\\
M_{31}\end{array}\right),\label{2.21}\end{equation} where
\[\tilde{M}=\left(\begin{array}{cc}(S_{1},S_{1})&M_{23}\\M_{32}&E\end{array}\right)-
\frac{1}{(S_{0},S_{0})}\left(\begin{array}{c}(S_{1},S_{0})\\M_{31}\end{array}\right)\left(\begin{array}{cc}
(S_{0},S_{1})&M_{13}\end{array}\right).\] Applying Lemma
\ref{matrix} in \eqref{2.21}, we get
\begin{equation}
I_{00}=1 +O\left(e^{-\frac{(\log
m)^{2}}{8}}\right).\label{2.17}\end{equation} In order to evaluate
the expansion of \eqref{2.3}, we are left to find
$\|S_{0}(x_{0})\|^{2}_{h_{m}}=\lambda_{0}^{2}$. From \eqref{2.20},
we have
\[\lambda_{0}^{2}=m(1+\frac{\rho}{2m})\left(1+O\big(e^{-(\log
m)^{2}}\big)\right).\]
 Therefore, the Tian-Yau-Zelditch expansion according to \eqref{2.3} on a Riemann surface with constant scalar curvature
 $\rho$ is
\begin{eqnarray*}
\lefteqn{I_{00}\lambda_{0}^{2}}\\&&=(1+O\left(e^{-\frac{(\log
m)^{2}}{8}}\right))m(1+\frac{\rho}{2m})\left(1+O\big(e^{-(\log
m)^{2}}\big)\right)\\&&=
m(1+\frac{\rho}{2m})+O\left(e^{-(\frac{(\log m)^{2}}{8})}\right)
\end{eqnarray*}
for
$m>\max\{e^{20\sqrt{5}}+2|\rho|,|\rho|^{4/3},\frac{1}{\delta},\sqrt{\frac{2}{|\rho|}}\}$.
\end{proof}

\end{document}